\documentclass[12pt]{amsart}
\usepackage{amsmath, amsthm, amsfonts, amssymb}

\newtheorem{theorem}{Theorem}[section]
\newtheorem{lemma}[theorem]{Lemma}

\theoremstyle{definition}

\theoremstyle{remark}

\numberwithin{equation}{section}

\newcommand{\mC}{\ensuremath{\mathbb{C}}}
\newcommand{\mD}{\ensuremath{\mathbb{D}}}

\newcommand{\mN}{\ensuremath{\mathbb{N}}}
\newcommand{\mR}{\ensuremath{\mathbb{R}}}

 \allowdisplaybreaks
 
\begin{document}

\title{Some Normality Criteria for Families of Holomorphic Functions of  Several Complex Variables}

\author[K. S. Charak]{Kuldeep Singh Charak}
\address{
\begin{tabular}{lll}
&Kuldeep Singh Charak\\
&Department of Mathematics\\
&University of Jammu\\
&Jammu-180 006\\ 
&India\\
\end{tabular}}
\email{kscharak7@rediffmail.com}

\author[R. Kumar]{Rahul Kumar}
\address{
\begin{tabular}{lll}
&Rahul Kumar\\
&Department of Mathematics\\
&University of Jammu\\
&Jammu-180 006\\
&India
\end{tabular}}
\email{rktkp5@gmail.com}
 
\begin{abstract}
We prove a Zalcman-Pang lemma in several complex variables and apply it to obtain several complex variables analogues of the known normality criteria like Lappan's five-point theorem and Schwick's theorem.
\end{abstract}

\renewcommand{\thefootnote}{\fnsymbol{footnote}}

\maketitle

\section{\textbf{Introduction}}
Let $D$ be a domain in $\mC^n$ and $\mathcal{F}$ be a family of holomorphic functions $f:D\rightarrow \mC.$ $\mathcal{F}$ is said to be normal in $D$ if every sequence in $ \mathcal{F}$ contains a subsequence that converges locally uniformly to a limit function which is either  holomorphic on $D$ or identically equal to $ \infty.$  $\mathcal{F}$ is said to be normal at a point $ z_{0} \in D $ if it is normal in some neighborhood of $z_{0}$ in $ D $.  As an attempt to obtain a natural extension of the theory of normal families of holomorphic functions of one complex variable (see \cite{Schiff, Zalcman}) to several complex variables, Dovbush\cite{Dov1} defined  the spherical derivative of a holomorphic function of several complex variables by using {\it Levi's form  } as follows:\\
For every $\psi\in \mathcal{C}^2(D),$ at each point $z$ of $D$ define a Hermitian form
\begin{equation}\label{eq1}
L_z(\psi,~v):= \sum\limits_{k,l=1}^{n} \frac{\partial^2 \psi}{\partial z_k\partial \bar{z_l}}(z)v_k\bar{v_l}
\end{equation} 
and is called the {\it Levi form } of the function $\psi$ at $z.$\\
For a holomorphic function $f$ defined on $D,$ define  
\begin{equation}\label{eq:2}
f^{\#}(z):= \sup\limits_{|v|=1}\sqrt{L_z(\log(1+|f|^2),~v)}.
\end{equation} 
Since $L_z(\log(1+|f|^2),~v)\geq 0,  \  f^{\#}(z)$ given by (\ref{eq:2}) is well defined  and for $n=1$ the formula $\eqref{eq:2}$ takes the form 
$$f^{\#}(z):= \frac{|f^\prime(z)|}{1+|f(z)|^2}$$
which is the spherical derivative on $\mC.$ Hence (\ref{eq:2}) gives the natural extension of the spherical derivative to $\mC^n.$\\
Also from \eqref{eq:2}, we find that
\begin{equation}\label{eq3}
f^{\#}(z) = \sup\limits_{|v|=1} \frac{|Df(z)v|}{1+|f(z)|^2}
\end{equation}
where $$D= (\frac{\partial}{\partial z_1},\frac{\partial}{\partial z_2}, \ldots ,\frac{\partial}{\partial z_n})$$

A well known powerful tool in the theory of normal families of holomorphic functions of one complex variables is the following lemma due to Zalcman\cite{Zalcman1}:

\medskip

{\bf Zalcman Lemma:} {\it A family $\mathcal{F}$ of holomorphic functions on the open unit disk $\mD$ is not normal in $\mD$ if and only if there exist a number $r: 0 < r < 1;$ points $z_n \in \{z:|z|<r\};$ functions $f_n\in \mathcal{F};$ and numbers $\rho_n \rightarrow 0$ such that 
$$g_n(\zeta)=f_n(z_n+\rho_n\zeta)\rightarrow g(\zeta), \ \mbox{ as } n\rightarrow \infty,$$
where $g$ is a nonconstant entire function satisfying $g^{\#}(\zeta)\leq g^{\#}(0)=1,$ for all $\zeta \in \mC.$}

\medskip

 Also, equally important is the following extension of the Zalcman Lemma due to Pang(\cite{Pang1}, Lemma $2$)( also\cite{Pang2}, Theorem $1$):

\medskip

{\bf Zalcman-Pang Lemma:} {\it Let $\mathcal{F}$ be a family of holomorphic functions on the open unit disk $\mD$ and $-1<\alpha<1$. Then $\mathcal{F}$ is not normal in $\mD$ if and only if there exist a number $r: 0 < r < 1;$ points $z_n \in \{z:|z|<r\};$ functions $f_n\in \mathcal{F};$ and numbers $\rho_n \rightarrow 0$ such that 
$$g_n(\zeta)=\rho_n^{-\alpha}f_n(z_n+\rho_n\zeta)\rightarrow g(\zeta), \ \mbox{ as } n\rightarrow \infty,$$
where $g$ is a nonconstant entire function satisfying $g^{\#}(\zeta)\leq g^{\#}(0)=1, \ \forall \zeta \in \mC.$}

\medskip

Dovbush\cite{Dov1} besides extending Marty's theorem\cite{Marty} extended  Zalcman Lemma to several complex variables as
\begin{theorem}(Zalcman Lemma in  $\mC^n$) Suppose that a family $\mathcal{F}$ of functions holomorphic on $D\subseteq\mathbb{C}^n$ is not normal at some point $w_0 \in D.$ Then there exist sequences $f_j \in \mathcal{F}, \ w_j\to \ w_0, \ \rho_j=1/f_j^{\#}(w_j)\to 0,$ such that the sequence $g_j(z)=f_j(w_j+\rho_j z) $ converges locally uniformly in $\mathbb{C}^n$ to a nonconstant entire function $g$ satisfying $g^{\#}(z)\leq g^{\#}(0)=1$ for all $z\in \mC^n.$
\label{ZLCN}
\end{theorem}
In this paper we give a several complex variables analogue of Zalcman-Pang Lemma, a generalization of Theorem \ref{ZLCN} and as applications, obtain several complex variables versions of Lappan's five-point theorem \cite{Lappan}, Schwick's theorem \cite{Schwick} and some other normality criteria.  

\section{\textbf{Main Results}}

\begin{theorem}(Zalcman-Pang Lemma in $\mC^n$)\label{thm:1}
Let $\mathcal{F}$ be a family of holomorphic functions on $D=\{z\in\mathbb{C}^n:|z|<1\}$ . If $\mathcal{F}$ is not normal on $D,$ then for all $\alpha :0\leq \alpha<1,$ there exist real number $r:0<r<1,$ and sequences $\{z_j\}\subseteq D: |z_j|<r,$  $\{f_j\}\subseteq\mathcal{F},$ and $ \{\rho_j\}\subset (0,\ 1]: \rho_j\to 0 $ such that
$$g_j({\zeta})= \frac{f_j(z_j+\rho_j\zeta)}{\rho_j^{\alpha}} $$ converges locally uniformly to a nonconstant entire function $g$ in $\mathbb{C}^n.$
\end{theorem}

\begin{theorem}\label{thm:2}
Let  $\mathcal{F}$ be a family of holomorphic functions on $D=\{z:|z|<1\}\subseteq\mathbb{C}^{n}$ and let  $\alpha$ and $\beta$ be real numbers such that $\alpha\geq 0$ and $\beta\geq \alpha+1.$ Then $\mathcal{F}$ is not normal on $D$ if and only if there exist real number $r:0<r<1$ and sequences $\{z_j\}\subseteq D: |z_j|<r,$  $\{f_j\}\subseteq\mathcal{F},$ and $\{\rho_j\}\subset (0, \ 1]: \rho_j\to 0 $ such that
$$g_j({\zeta})= \rho_j^{-\alpha}f_j(z_j+\rho_j^{\beta}\zeta) $$ 
converges locally uniformly to a nonconstant entire function $g$ in $\mathbb{C}^n.$ 
\end{theorem}

For $\alpha=0$ and $\beta=1$, Theorem\ref{thm:2} reduces to  Theorem\ref{ZLCN}.

\medskip

By Theorem\ref{thm:2}, we extend Lappan's five-point theorem\cite{Lappan}((also see,  Hinkkanen\cite{Hink}) to several complex variables as
\begin{theorem}\label{thm2} A family $\mathcal{F}$ of holomorphic functions on a domain $D\subseteq \mathbb{C}^{n}$ is normal on $D$ if and only if there exists a set $E$ containing at least three points such that for each compact subset $K\subset D,$ there exists a positive constant $M(K)$ for which 
\begin{equation}\label{eq3}
f^{\sharp}(z)\leq M(K) \mbox{ whenever }  f(z)\in E, ~ z\in K, ~ f\in\mathcal{F}.
\end{equation}
\end{theorem}

Schwick \cite{Schwick} sharpened Royden's theorem \cite{Royden} as: {\it Let $\mathcal{F}$ be a family of meromorphic functions on a domain $D$ with the property that for each compact set $K \subset D$ there is a function $h_{K} : [0,\infty] \rightarrow [0,\infty]$, which is finite somewhere on $(0,\infty)$, such that 
\begin{equation}
 |f^{'}(z)| \leq h_{K}(|f(z)|), \ \mbox{ for all } f \in \mathcal{F}, z \in K.
\label{alpha}
\end{equation}
 Then $\mathcal{F}$ is normal on $D.$} Actually, Schwick's theorem requires (\ref{alpha}) to be satisfied by $f$ at least on a circle. Wang\cite{Wang}, by applying Zalcman's lemma, obtained the following more sharpened version of Schwick's theorem wherein (\ref{alpha}) is required to be satisfied by $f$ at least for five points:
\begin{theorem} Let $\mathcal{F}$ be a family of meromorphic functions on a domain $D\subset \mC$ with the property that for each compact set $K \subset D$ there is a function $h_{K} : \overline{\mC} \rightarrow [0,\infty]$, which is finite for at least five points on $\overline{\mC}$, such that 
\begin{equation}
 |f^{'}(z)| \leq h_{K}(f(z)), \ \mbox{ for all } f \in \mathcal{F}, z \in K.
\label{beta}
\end{equation}
 Then $\mathcal{F}$ is normal on $D.$\\
Moreover, a family $\mathcal{F}$ of holomorphic functions is normal on $D,$ if (\ref{beta}) is satisfied and the function $h_K$ is finite for at least three points on $\mC.$
\label{Wang}
\end{theorem} 
By using Theorem \ref{thm2}, we obtain a several complex variables analogue of Theorem \ref{Wang}:

\begin{theorem}\label{thm3} Let $\mathcal{F}$ be a family of holomorphic functions on a domain $D\subseteq\mathbb{C}^{n}$ with the property that for each compact subset $K\subset D$ there is a function $h_{K}:\overline{\mC}\longrightarrow[0,~\infty]$, which is finite for at least three points on $\mC$ such that $\left|Df(z)\right|\leq h_{K}(f(z))$ for all $f\in \mathcal{F}$ and $z\in K.$ Then $\mathcal{F}$ is normal on $D.$  
\end{theorem}
 
Further, we obtain a several complex variables version of a normality criterion due to  Tan and Thin(\cite{Tan}, Theorem $1$, page $48$).
For the sake of convenience,  we shall use the following notations:

 $$f_{z_j}=\frac{\partial f}{\partial z_j} \ \  \mbox{ and } \ \ f_{z_kz_j}=\frac{\partial^2 f }{\partial z_j \partial z_k}.$$

\begin{theorem}\label{thm4}
Let $\mathcal{F}$ be a family of holomorphic functions on a domain $D\subseteq\mathbb{C}^{n}.$ Assume that for each compact subset $K\subset D,$ there exist a  set $E=E(K)\subset\mathbb{C}$ consisting of two distinct points and a positive constant $M=M(K)$ such that 
$$f^{\sharp}(z)\leq M  \mbox{ and } (f_{z_k})^{\sharp}(z)\leq M,\mbox{ whenever } z\in K,~f(z)\in E,~ k=1, 2, \ldots, n.$$
Then $\mathcal{F}$ is normal on $D.$
\end{theorem}

Finally, we obtain another version of a result due to Cao and Liu(\cite{Cao}, Theorem $1.8(i)$ page $1395$):
\begin{theorem}\label{thm5}
Let $\mathcal{F}$ be a family of holomorphic functions in a domain $D= \{z\in \mathbb{C}^n:|z|<1\}$ and $s>0$ be any real number. If
$$\mathcal{G}= \{\frac{|Df(z)|}{1+|f(z)|^s}: f\in \mathcal{F}\}$$
is locally uniformly bounded in $D,$ then $\mathcal{F}$ is normal on $D.$
\end{theorem}

\section{\textbf{Proofs of Main Results}}

Let $f$ be a meromorphic function in $\mC$ and $ a\in \overline{\mC}.$ Then $a$ is called totally ramified value of $f$ if $f-a$ has no simple zeros. Following result known as {\it Nevanlinna's Theorem } (see \cite{Bergweiler}) plays a crucial role in our proofs:
\begin{theorem}\label{thm1}
Let $f$ be a non-constant meromorphic function $a_1, a_2, \ldots, a_q\in \overline{\mC} \mbox{ and } m_1, m_2, \dots,  m_q \in \mN.$ Suppose that all $a_j$-points of $f$ have multiplicity  at least $m_j,  \mbox{ for }~ j=1, 2, \ldots, q.$ Then 
$$\sum\limits_{j=1}^{q}(1-\frac{1}{m_j}) \leq 2.$$
\end{theorem}
If $f$ does not assume  the value $a_j$ at all, then  we take  $m_j=\infty.$ From Theorem \ref{thm1}, it follows that if $f$ is entire function and $a_1, ~a_2 \in \mC$ are distinct such that all $a_j-$ points of $f$ have multiplicity at least $3,$ then $f$ is constant. Also, it follows that if $a_1,\ a_2,\ a_3 \in \mC$ are distinct such that all $a_j-$ points of $f$ have multiplicity at least $2,$ then $f$ is constant. Thus, a non-constant entire function can not have more than two totally ramified values. 

\medskip

 For the proof of Theorem \ref{thm:1} we need the following lemma:

\begin{lemma}\label{lemma01}
Let $f$ be a holomorphic function in $D=\{z\in\mathbb{C}^n:|z|<1\}$ and let $-1<\alpha<1.$ Let $\Omega:=\{z:|z|<r<1\}\times (0, \ 1]$ and $F:\Omega\rightarrow \mR$ be defined as
$$F(z,t)=\frac{(r-|z|)^{1+\alpha}t^{1+\alpha}(1+|f(z)|^2)f^{\sharp}(z)}{(r-|z|)^{2\alpha}t^{2\alpha}+|f(z)|^2}.$$
If $F(z,1)>1$ for some $z\in \{z:|z|<r<1\},$ then there exist $z_0\in \{z:|z|<r<1\}$ and $t_0 \in (0, \ 1)$ such that 
$$\sup_{|z|< r}F(z,t_0)=F(z_0,t_0)=1.$$
\end{lemma}

A small variation in Lemma\ref{lemma01} yields:

\begin{lemma}\label{lemma1}
Let $f$ be a holomorphic function in $D=\{z\in\mathbb{C}^n:|z|<1\}$ and let $0\leq \alpha<\beta.$ Let $\Omega:=\{z:|z|<r<1\}\times (0, \ 1]$ and $F:\Omega\rightarrow \mR$ be defined as
$$F(z,t)=\frac{(r-|z|)^{\beta+\alpha}t^{\beta+\alpha}(1+|f(z)|^2)f^{\sharp}(z)}{(r-|z|)^{2\alpha}t^{2\alpha}+|f(z)|^2}.$$
If $F(z,1)>1$ for some $z\in \{z:|z|<r<1\},$ then there exist $z_0\in \{z:|z|<r<1\}$ and $t_0 \in (0, \ 1)$ such that 
$$\sup_{|z|< r}F(z,t_0)=F(z_0,t_0)=1.$$
\end{lemma}

\textbf{Proof of Lemma \ref{lemma01}:} First, we show that 
\begin{equation}
\lim\limits_{(r-|z|)t\to 0}F(z,t)=0.
\label{00}
\end{equation}

Since $F$ is continuous on $\Omega$, we shall prove $(\ref{00})$ for $(r-|z|)t\to 0$ through an arbitrary sequence 
$x_j=(r-|z_j|)t_j\to 0 ~\mbox{as}~ j\to\infty$ where $z_j\in\{z:|z|<r\}, ~t_j\in(0,~1).$ Put $\lim_{j\to\infty}z_j=w_0.$ Then $|w_0|\leq r.$\\
If $f(w_0)\neq 0,$ then for $-1<\alpha,$ we have
\begin{eqnarray*}
0 &\leq& \lim\limits_{j\to\infty}F(z_j,t_j)\\
  &\leq& \lim\limits_{j\to\infty}\frac{x_j^{1+\alpha}(1+|f(z_j)|^2)f^{\sharp}(z_j)}{|f(z_j)|^2}\\
	&=& 0
\end{eqnarray*}

If $f(w_0)= 0,$ then for $\alpha<1,$ we have
\begin{eqnarray*}
0 &\leq& \lim\limits_{j\to\infty}F(z_j,t_j)\\
  &\leq& \lim\limits_{j\to\infty}x_j^{1-\alpha}(1+|f(z_j)|^2)f^{\sharp}(z_j)\\
	&=& 0
\end{eqnarray*}

Hence $(\ref{00})$ holds.\\
Let 
$$U:= \{(z,t)\in \Omega: F(z,t)>1\}.$$
Since $F(z,1)>1$ for some $z=z^*\in \{z:|z|<r<1\}$, $U\neq \emptyset.$ Clearly, $t_0:= \inf\{t: (z,t)\in U\}\neq 0.$ Also, $t_0\neq 1$  since otherwise there exists a sequence $\{t_j\} ~(<1)$ such that $t_j\to t_0$ as $j\to\infty$ and $F(z^*,t_j)\leq 1.$ This implies that 
$$\lim_{j\to\infty}F(z^*,t_j)=F(z^*,1)\leq 1,$$ 
which contradicts that $(z^*,1)\in U.$ Hence $0<t_0<1.$

\smallskip

Now we take $z_0\in\{z:|z|\leq r\} $ such that  
$$\sup_{|z|\leq r} F(z,t_0)= F(z_0, t_0). $$ 
 To complete the proof we shall show that $F(z_0,t_0)=1.$ Suppose this is not true. Then we have the following two cases:

\smallskip

{\it Case $1:$} When $F(z_0, t_0)<1.$ In this case there exists a sequence $(z_j,t_j)\in U$ such that $~t_j\to t_0.$ Let $z_j\to w_1.$ Then $|w_1|\leq r.$ Since $F(w_1,t_0)\leq F(z_0,t_0)<1,$ by continuity of $F$ it follows that for  
for sufficiently large $j,$  $F(z_j,t_j)<1,$ a contradiction. 

\smallskip

{\it Case $2:$} When $F(z_0, t_0)>1.$ Since $F(z_0,0)=0,$ by continuity of $F$ with respect to $t,$ there exists $t_1:0<t_1<t_0$ such that
$$F(z_0,t_1)= 1+\frac{F(z_0,t_0)-1}{2}$$
which contradicts the definition of $t_0.$ $\Box$

\medskip

\textbf{Proof of Theorem \ref{thm:1}:} Without loss of generality, we may assume that $D=\{z:|z|<1\}$ and let $\mathcal{F}$ be not normal at $z_0 = 0.$ Then by several complex variables analogue of Marty's theorem(see \cite{Dov1}, Theorem$2.1$), there exist $r_0:0<r_0<1, ~\{w_j\}\subset \{z:|z|<r_0\},$ and $\{f_j\}\subseteq\mathcal{F}$ such that 
$$\lim\limits_{j\to\infty} f_j^\sharp(w_j)=\infty.$$
Choose $r$ such that $0<r_0<r<1$ and corresponding to each $f_j \in \mathcal{F}$ define $F_j:\{z:|z|<r\}\times (0,~1]\rightarrow \mR$ as 
$$F_j(z,t)= \frac{(r-|z|)^{1+\alpha}t^{1+\alpha}(1+|f_j(z)|^2)f_j^\sharp(z)}{(r-|z|)^{2\alpha}t^{2\alpha}+|f_j(z)|^2}.$$ 
Then 
\begin{eqnarray*}
F_j(w_j,1) &=& \frac{(r-|w_j|)^{1+\alpha}(1+|f_j(w_j)|^2)f_j^\sharp(w_j)}{(r-|w_j|)^{2\alpha}+|f_j(w_j)|^2}\\
            &=& \frac{(r-|w_j|)^{1-\alpha}(1+|f_j(w_j)|^2)f_j^\sharp(w_j)}{1+\frac{|f_j(w_j)|^2}{(r-|w_j|)^{2\alpha}}}\\
						&>& \frac{(r-r_0)^{1-\alpha}(1+|f_j(w_j)|^2)f_j^\sharp(w_j)}{1+\frac{|f_j(w_j)|^2}{(r-r_0)^{2\alpha}}} \to \infty ~\mbox{as}~ j\to\infty
\end{eqnarray*}
Thus for sufficiently large  $j,$  $F_j(w_j,1)>1.$ and hence by Lemma{\ref{lemma01}}, there exist $z_j\in \{z:|z|<r\} ~\mbox{and}~ t_j\in (0,~1)$ such that 
$$\sup_{|z|< r}F_j(z,t_j)= F_j(z_j,t_j)=1.$$
Thus, for sufficiently large $j$, we have 
\begin{eqnarray*}
1 &=& F_j(z_j, t_j)\\
  &\geq& F_j(w_j, t_j)\\
	&=&  \frac{(r-|w_j|)^{1+\alpha}t_j^{1+\alpha}(1+|f_j(w_j)|^2)f_j^\sharp(w_j)}{(r-|w_j|)^{2\alpha}t_j^{2\alpha}+|f_j(w_j)|^2}\\
	&\geq& \frac{t_j^{1+\alpha}(r-|w_j|)^{1+\alpha}(1+|f_j(w_j)|^2)f_j^\sharp(w_j)}{(r-|w_j|)^{2\alpha}+|f_j(w_j)|^2}\\
	&=& t_j^{1+\alpha}F_j(w_j,1)
\end{eqnarray*}
 which implies that $\lim\limits_{j\to\infty}t_j=0.$ 
Let $\rho_j=(r-|z_j|)t_j\to 0.$ Then
$$\lim\limits_{j\to\infty}\frac{\rho_j}{r-|z_j|}=0.$$
Thus the function 
$$g_j(\zeta):= \frac{f_j(z_j+\rho_j\zeta)}{\rho_j^{\alpha}}$$
is defined for 
$$|\zeta|<R_j=\frac{r-|z_j|}{\rho_j}\to\infty.$$
Now 
\begin{eqnarray}\label{*}
\sup_{|v|=1}\frac{|Dg_j(\zeta)v|}{1+|g_j(\zeta)|^2} &=& \sup_{|v|=1}\frac{\rho_j^{1-\alpha}|Df_j(z_j+\rho_j\zeta)v|}{1+\frac{|f_j(z_j+\rho_j\zeta)|^2}{\rho_j^{2\alpha}}}\nonumber\\
                                                   &=& \sup_{|v|=1}\frac{\rho_j^{1+\alpha}|Df_j(z_j+\rho_j\zeta)v|}{\rho_j^{2\alpha}+|f_j(z_j+\rho_j\zeta)|^2} 
\end{eqnarray}
Since
 $$\frac{r-|z_j|}{r-|z_j+\rho_j\zeta|}\to 1,$$
there exists $\epsilon_j \to 0$ such that 
$$\rho_j^{1+\alpha}\leq (1+\epsilon_j)^{1+\alpha}(r-|z_j+\rho_j\zeta|)^{1+\alpha}t_j^{1+\alpha}, $$
and
$$\rho_j^{2\alpha}\geq(1-\epsilon_j)^{2\alpha}(r-|z_j+\rho_j\zeta|)^{2\alpha}t_j^{2\alpha}.$$
Thus from $(\ref{*}),$ we get
\begin{eqnarray*}
\sup_{|v|=1}\frac{|Dg_j(\zeta)v|}{1+|g_j(\zeta)|^2} &\leq& \sup_{|v|=1}\frac{(1+\epsilon_j)^{1+\alpha}(r-|z_j+\rho_j\zeta|)^{\alpha+1}t_j^{1+\alpha}|Df_j(z_j+\rho_j\zeta)v|}{(1-\epsilon_j)^{2\alpha}(r-|z_j+\rho_j\zeta|)^{2\alpha}t_j^{2\alpha}+|f_j(z_j+\rho_j\zeta)|^2}\\
&=& \frac{(1+\epsilon_j)^{1+\alpha}(r-|z_j+\rho_j\zeta|)^{1+\alpha}t_j^{1+\alpha}(1+|f_j(z_j+\rho_j\zeta)|^2)f_j^{\sharp}(z_j+\rho_j\zeta)}{(1-\epsilon_j)^{2\alpha}(r-|z_j+\rho_j\zeta|)^{2\alpha}t_j^{2\alpha}+|f_j(z_j+\rho_j\zeta)|^2}\nonumber\\
&\leq& \frac{(1+\epsilon_j)^{1+\alpha}}{(1-\epsilon_j)^{2\alpha}}
\end{eqnarray*}
That is, 
 $$g_j^{\sharp}(\zeta)\leq \frac{(1+\epsilon_j)^{1+\alpha}}{(1-\epsilon_j)^{2\alpha}}$$

and hence by Marty's theorem $\{g_j\}$ is normal in $\mathbb{C}^n.$ Without loss of generality we may assume that $\{g_j\}$ converges locally uniformly to a holomorphic function $g$ or $\infty$ in $\mathbb{C}^n.$
Since  
\begin{eqnarray*}
g_j^{\sharp}(0) &=& \sup_{|v|=1}\frac{|Dg_j(0)v|}{1+|g_j(0)|^2}\\
                &=& \sup_{|v|=1} \frac{\rho_j^{1+\alpha}|Df_j(z_j)v|}{\rho_j^{2\alpha}+|f_j(z_j)|^2}\\
								&=& \frac{(r-|z_j|)^{1+\alpha}t_j^{1+\alpha}f_j^{\sharp}(z_j)(1+|f_j(z_j)|^2)}{(r-|z_j|)^{2\alpha}t_j^{2\alpha}+|f_j(z_j)|^2}\\
						    &=& F_j(z_j,t_j) = 1,
\end{eqnarray*}
it follows that $g(\zeta)$ is a nonconstant entire function in $\mathbb{C}^n.$ $\Box$

\medskip

\textbf{Proof of Theorem \ref{thm:2}:} Let  $\mathcal{F}$ be a family of holomorphic functions on $D=\{z:|z|<1\}\subseteq\mathbb{C}^{n}$ and let  $\alpha$ and $\beta$ be real numbers such that $\alpha\geq 0$ and $\beta\geq \alpha+1.$ Further, suppose that there exist real number $r:0<r<1$ and sequences $\{z_j\}\subseteq D: |z_j|<r,$  $\{f_j\}\subseteq\mathcal{F},$ and $\{\rho_j\}\subset (0, \ 1]: \rho_j\to 0 $ such that
$$g_j({\zeta})= \rho_j^{-\alpha}f_j(z_j+\rho_j^{\beta}\zeta) $$ 
converges locally uniformly to a nonconstant entire function $g$ in $\mathbb{C}^n.$ Then There is some $\zeta_0\in\mathbb{C}^n$ such that $g^{\sharp}(\zeta_0)>0.$ Suppose $z_j\to z_0$ as $j\to\infty.$ Then $|z_0|\leq r.$ Since  
$$|g_{j_{z_1}}(\zeta_0).v_1+ \ldots +g_{j_{z_n}}(\zeta_0).v_n|= \rho_j^{\beta-\alpha}|f_{j_{z_1}}(z_j+\rho_j^{\beta}\zeta_0).v_1+ \ldots +f_{j_{z_n}}(z_j+\rho_j^{m}\zeta_0).v_n|, $$
it follows that 
\begin{eqnarray*}
f_j^{\sharp}(z_j+\rho_j^{\beta}\zeta_0) &=& \sup_{|v|=1}\frac{|f_{j_{z_1}}(z_j+\rho_j^{\beta}\zeta_0).v_1+ \ldots +f_{j_{z_n}}(z_j+\rho_j^{\beta}\zeta_0).v_n|}{1+|f_j(z_j+\rho_j^{\beta}\zeta_0)|^2}\\
                                      &=& \sup_{|v|=1}\frac{\rho_j^{\alpha-\beta}|g_{j_{z_1}}(\zeta_0).v_1+ \ldots +g_{j_{z_n}}(\zeta_0).v_n|}{1+\rho_j^{2\alpha}|g_j(\zeta_0)|^2}\\
																			&\geq& \sup_{|v|=1}\frac{\rho_j^{\alpha-\beta}|g_{j_{z_1}}(\zeta_0).v_1+ \ldots +g_{j_{z_n}}(\zeta_0).v_n|}{1+|g_j(\zeta_0)|^2}\\
																			&=& \rho_j^{\alpha-\beta}g_j^{\sharp}(\zeta_0)\to\infty ~\mbox{as}~ j \to\infty
\end{eqnarray*}
and so by Marty's theorem $\mathcal{F}$ is not normal at $z_0$ and hence $\mathcal{F}$ is not normal on $D.$

\smallskip

Conversely, suppose that $\mathcal{F}$ is not normal at $z_0 = 0.$ Then by Marty's theorem, there exist $0<r^*<1, ~|z_j^*|<r^*, ~\{f_j\}\subseteq\mathcal{F}$ such that 
$$\lim\limits_{j\to\infty} f_j^\sharp(z_j^*)=\infty.$$
Choose $r$ such that $0<r^*<r<1$ and corresponding to each $f_j\in\mathcal{F}$ define 
$$F_j(z,t):= \frac{(r-|z|)^{\beta+\alpha}t^{\beta+\alpha}(1+|f_j(z)|^2)f_j^\sharp(z)}{(r-|z|)^{2\alpha}t^{2\alpha}+|f_j(z)|^2},$$ 
where $|z|<r, ~0<t\leq 1.$  Then
\begin{eqnarray*}
F_j(z_j^*,1) &=& \frac{(r-|z_j^*|)^{\beta+\alpha}(1+|f_j(z_j^*)|^2)f_j^\sharp(z_j^*)}{(r-|z_j^*|)^{2\alpha}+|f_j(z_j^*)|^2}\\
            &=& \frac{(r-|z_j^*|)^{\beta-\alpha}(1+|f_j(z_j^*)|^2)f_j^\sharp(z_j^*)}{1+\frac{|f_j(z_j^*)|^2}{(r-|z_j^*|)^{2\alpha}}}\\
						&>& \frac{(r-r^*)^{\beta-\alpha}(1+|f_j(z_j^*)|^2)f_j^\sharp(z_j^*)}{1+\frac{|f_j(z_j^*)|^2}{(r-r^*)^{2\alpha}}} \to \infty ~\mbox{as}~ j\to\infty.
\end{eqnarray*}
Thus for large $j$, we have
$$F_j(z_j^*,1)>1$$
and therefore, by Lemma {\ref{lemma1}}, there exist ${z_j} ~\mbox{and}~ {t_j}$ satisfying $|z_j|< r, ~0<t_j<1$  such that 
$$\sup_{|z|< r}F_j(z,t_j)= F_j(z_j,t_j)=1.$$
Now rest of the proof goes on the same lines as that of the proof of Theorem\ref{thm:1}. $\Box$

\medskip
 
\textbf{Proof of Theorem \ref{thm2}:} By Marty's theorem in $\mC^n$ (see \cite{Dov1}, Theorem${2.1}$) we find that  \eqref{eq3} is necessary with $E= \mC.$ To prove the sufficiency, suppose \eqref{eq3} holds but $\mathcal{F}$ is not normal. Then by Theorem\ref{ZLCN}  there exist sequences $\{f_j\}\subset \mathcal{F};$ ~$\{w_j\}\subset D : w_j\rightarrow w_0$  and $\{\rho_j\}\subset (0,~1):\rho_j \rightarrow 0,$ such that the sequence $\{g_j\}$ defined as $g_j(\zeta )=f_j(w_j+\rho_j\zeta)$ converges  locally uniformly on $\mC^n$ to a nonconstant entire function $g.$ Let $K$ be a compact set containing $w_0$ and suppose $g(\zeta_0)\in E$. By Hurwitz's theorem, there exists $\zeta_j\rightarrow\zeta_0$ such that $f_j(w_j+\rho_j\zeta_j)=g_j(\zeta_j )=g(\zeta_0)~\mbox{for large}~j.$ 
Since $f^{\sharp}_j(w_j+\rho_j\zeta_j)\leq M$ for $j$ sufficiently large, we have
$$g^{\sharp}(\zeta_0)= \lim\limits_{j\to\infty} g_j^{\sharp}(\zeta_j)= \lim\limits_{j\to\infty} \rho_j f^{\sharp}_j(w_j+\rho _j\zeta_j)\leq \lim\limits_{j\to\infty}\rho_j M= 0.$$
Thus $g^{\sharp}(\zeta_0)=0$ whenever $g(\zeta_0)\in E$ implying that
$$\sup\limits_{|v|=1}\left[\frac{\left|g_{z_1}(\zeta_0)\cdot v_1+g_{z_2}(\zeta_0)\cdot v_2+ \ldots +g_{z_n}(\zeta_0)\cdot v_n\right|}{1+ |g(\zeta_0)|^2}\right] =0 $$
whenever $g(\zeta_0)\in E$ which further implies that  
$$g_{z_1}(\zeta_0)\cdot v_1 + g_{z_2}(\zeta_0)\cdot v_2 + \ldots + g_{z_n}(\zeta_0)\cdot v_n=0$$
whenever $g(\zeta_0)\in E,$ for all $(v_1,~v_2, \ldots,~v_n)$ such that 
$$\sqrt{|v_1|^2+|v_2|^2+ \ldots +|v_n|^2}=1.$$
Taking  $v_k=1$ and $v_m=0$ for all $m\neq k.$ Then $g_{z_k}(\zeta_0)=0 \mbox{ whenever } g(\zeta_0)\in E.$
Now, let $w= (a_1,~a_2,\ldots,~a_n), ~ w'=(b_1,~b_2,\ldots,~b_n)\in\mC^n$ and define 
$$h_j(z_j):= g(b_1,\ldots, ~b_{j-1},~z_j,~a_{j+1}, \ldots ,~a_n), \ \ j=1,2, \ldots, n.$$

Suppose $h_j(a)\in E.$ Then $ g(b_1,\ldots,~b_{j-1},~a,~a_{j+1}, \ldots,~a_n)\in E$ and hence 
$$g_{z_j}(b_1,\ldots,b_{j-1},~a,~a_{j+1},\ldots,~a_n)=0.$$
That is, 
$$\frac{dh_j}{dz_j}(a)=0, \ j=1,2, \ldots, n.$$ 
This, by Theorem \ref{thm1}, implies that each $h_j(z_j)$ is constant. Thus for $j=1,2,\ldots, n, $ we have
$$h_j(z_j)= g(b_1,\ldots,b_{j-1},~z_j,~a_{j+1},\ldots,~a_n)= \mbox{ a constant }$$ 
which implies that $g(w)=g(w')$ for all  $w,~w'\in\mC^n$ showing that $g$ is constant, a contradiction. $\Box$

\medskip

\textbf{Proof of Theorem \ref{thm3}:} Let $K$ be a compact subset of $D$ and let $\zeta_1, \zeta_2, \zeta_3 \in \mC$ be such that $h_{K}(\zeta_j)< \infty.$ Put 
$$E(K)=\{\zeta_1,\ \zeta_2, \ \zeta_3\} \mbox{ and } M(K)= \max_{\zeta\in E}|h_K(\zeta)|.$$
Then, for each $f\in \mathcal{F},$ we have
\begin{eqnarray*}
f^{\sharp}(z)&=&\sup\limits_{\left|v\right|=1}\frac{\left|Df(z)v\right|}{1+\left|f(z)\right|^{2}}\\
             &\leq& \sup\limits_{\left|v\right|=1}\left|Df(z)\right|\left|v\right|\\
             &=&\left|Df(z)\right| \leq h_{K}(f(z))\leq M(K) 
\end{eqnarray*}
whenever $ z\in K$ and $f(z)\in E(K).$\\
By Theorem \ref{thm2}  the normality of $\mathcal{F}$ on $D$ follows. $\Box$

\medskip

\textbf{Proof of Theorem \ref{thm4}:} Suppose $\mathcal{F}$ is not normal. Then, by Theorem \ref{ZLCN}, there exist sequences $f_{j}\in \mathcal{F},~w_{j}\to w_{0},~\rho_{j} \to 0,$ such that the sequence $g_{j}(\zeta)=f_{j}(w_{j}+\rho_{j}\zeta)$ converges locally uniformly in $\mC^{n}$ to a non-constant entire function $g.$ Let $K$ be a compact set containing $w_{0}$. Then there exists a set $E$ containing two points and $M>0$ such that $f^{\sharp}(z)\leq M, ~(f_{z_{k}})^{\sharp}(z)\leq M$  whenever $z\in K,~f(z)\in E.$ Let $w =(a_{1},a_{2},\ldots, a_{n}),~w' =(b_{1},b_{2}, \ldots,b_{n}) \in \mC^{n}$ and define 
$$h_{i}(z_{i}):=g(b_{1}, \ldots, b_{i-1},z_{i},a_{i+1}, \ldots,a_{n}), ~i= 1,2, \ldots, n.$$
First, we shall show that for any $a\in E$, all zeros of $h_{i}(z_{i})-a$ have multiplicity at least $3.$ Let $c$ be zero of $h_{i}(z_{i})-a$. Then $\zeta_{0}=(b_{1},\ldots, b_{i-1},c,a_{i+1}, \ldots, a_n)$ is a zero of $g(z)-a.$ By Hurwitz's theorem, there exists a sequence $\zeta_{j}\to \zeta_{0}$ such that $f_{j}(w_{j}+\rho_{j}\zeta_{j}) \to a$ and therefore, $w_{j}+\rho_{j}\zeta_{j}\in K$ and $f_{j}(w_{j}+\rho_{j}\zeta_{j})\in E \mbox{ for large } j. $ Hence 
$$f_{j}^{\sharp}(w_{j}+\rho_{j}\zeta_{j}) \leq M, ~(f_{j_{z_{k}}})^{\sharp}(w_{j}+\rho_{j}\zeta_{j})\leq M,~k=1,2, \ldots,n.$$
Now, 
\begin{eqnarray*}
g_{j}^{\sharp}(\zeta_{j}) &=& \sup_{|v|=1} \frac{|{g_j}_{z_1}(\zeta_{j}).v_{1} + \ldots +{g_j}_{z_n}(\zeta_{j}).v_{n}|}{1+|g_{j}(\zeta_{j})|^{2}} \\
                    &=& \sup_{|v|=1} \frac{\rho_{j}|f_{j_{z_{1}}}(w_{j}+\rho_{j}\zeta_{j}).v_{1}+ \ldots +f_{j_{z_{n}}}(w_{j}+\rho_{j}\zeta_{j}).v_{n}|}{1+|f_{j}(w_{j}+\rho_{j}\zeta_{j})|^{2}} \\
										&=& \rho_{j}f_{j}^{\sharp}(w_{j}+\rho_{j}\zeta_{j}) \\
										&\leq& \rho_{j}M \to 0 ~as ~j \to \infty.
\end{eqnarray*}										
Thus  $g^{\sharp}(\zeta_{0}) = 0 $ which implies that
$$\sup_{|v|=1} \frac{|g_{z_1}(\zeta_{0}).v_{1} + \ldots +g_{z_i}(\zeta_{0}).v_{i}+ \ldots +g_{z_n}(\zeta_{0}).v_{n}|}{1+|g(\zeta_{0})|^{2}} = 0$$
and hence $g_{z_i}(\zeta_{0}) = 0$ for each $i=1,2,\ldots, n$. That is, $ g_{z_i}(b_{1}, \ldots, b_{i-1},c,a_{i+1}, \ldots ,a_{n}) = 0$ implying that 
$$ \frac{dh_i}{dz_i}(c)= 0, \ i=1,2, \ldots, n.$$ 
This shows that $c$ is an $a-$point of $h_{i}$ with multiplicity at least $2.$\\
Next,
\begin{eqnarray*}
(g_{j_{z_i}})^{\sharp}(\zeta_{j}) &=&  \sup_{|v|=1} \frac{|g_{j_{z_iz_1}}(\zeta_{j}).v_{1} + \ldots +g_{j_{z_iz_n}}(\zeta_{j}).v_{n}|}{1+|g_{j_{z_{i}}}(\zeta_{j})|^{2}} \\
                               		&=& \sup_{|v|=1} \frac{\rho_{j}^{2}|f_{j_{z_{i}z_{1}}}(w_{j}+\rho_{j}\zeta_{j}).v_{1}+ \ldots +f_{j_{z_{i}z_{n}}}(w_{j}+\rho_{j}\zeta_{j}).v_{n}|}{1+\rho_{j^{2}}|f_{j_{z_{i}}}(w_{j}+\rho_{j}\zeta_{j})|^{2}} \\
															&=& \sup_{|v|=1} \frac{\rho_{j}^{2}(f_{j_{z_{i}}})^{\sharp}(w_{j}+\rho_{j}\zeta_{j})[1+|f_{j_{z_{i}}}(w_{j}+\rho_{j}\zeta_{j})|^{2}]}{1+\rho_{j}^{2}|f_{j_{z_{i}}}(w_{j}+\rho_{j}\zeta_{j})|^{2}} 
\end{eqnarray*}
Since\\
$\sup_{|v|=1} |f_{j_{z_{1}}}(w_{j}+\rho_{j}\zeta_{j}).v_{1}+ \ldots +f_{j_{z_{n}}}(w_{j}+\rho_{j}\zeta_{j}).v_{n}| = f_{j}^{\sharp}(w_{j}+\rho_{j}\zeta_{j})[1+|f_{j}(w_{j}+\rho_{j}\zeta_{j})|^{2}], $
therefore, $|f_{j_{z_{i}}}(w_{j}+\rho_{j}\zeta_{j})|	< M[1+\max_{d\in E}|d|^{2}].$ Thus,
\begin{eqnarray*}
(g_{j_{z_{i}}})^{\sharp}(\zeta_{j}) &\leq& \frac{\rho_{j}^{2}.M[1+\{M(1+\max_{d \in E}|d|^{2})\}^{2}]}{1+\rho_{j}^{2}|f_{j_{z_{i}}}(w_{j}+\rho_{j}\zeta_{j})|^{2}}\\
                               &\leq& M[1+\{M(1+\max\limits_{d \in E}|d|^{2})\}^{2}]\rho_{j}^{2}\\
															&\to & 0 ~as ~j \to\infty	
\end{eqnarray*}
and hence $(g_{z_{i}})^{\sharp}(\zeta_{0})=0.$
That is,
$$\sup_{|v|=1} |g_{z_iz_1}(\zeta_{0}).v_{1}+ \ldots +g_{z_iz_n}(\zeta_{0}).v_{n}| = 0.$$
That is, $g_{z_iz_i}(\zeta_{0}) = 0 $ implying that $g_{z_iz_i}(b_{1},\ldots,b_{i-1},c,a_{i+1},\ldots, a_{n})= 0 $. Hence
$$\frac{d^{2}}{dz_{i}^{2}}h_{i}(c)= 0 $$
 showing that	 $c$ is an $a-$point of $h_{i}(z_{i})$ with multiplicity at least $3$. Now by Theorem\ref{thm1}, we conclude that each $h_{i}(z_{i})$ is constant and hence $g$ is constant, a contradiction. $\Box$

\medskip

\textbf{Proof of Theorem \ref{thm5}:}  Suppose $\mathcal{F}$ is not normal. Then, by Theorem \ref{ZLCN}, there exist sequences $f_{j}\in \mathcal{F},~w_{j}\to w_{0},~\rho_{j} \to 0,$ such that the sequence $g_{j}(\zeta)=f_{j}(w_{j}+\rho_{j}\zeta)$ converges locally uniformly in $\mC^{n}$ to a non-constant entire function $g.$ \\
Let $K$ be compact set containing $w_0.$ Since $\mathcal{G}$ is locally uniformly bounded in $D.$ So there exist some constant $M(K)>0$ such that 
$$\frac{|Df(z)|}{1+|f(z)|^s}\leq M, ~z\in K, ~f\in\mathcal{F}.$$
Now \begin{eqnarray*}
     |Dg_j(\zeta)| &=& \rho_j|Df_j(w_j+\rho_j\zeta)|\\
		               &\leq& \rho_j.M(1+|f_j(w_j+\rho_j\zeta)|^s)\\
									 &=& \rho_j. M(1+|g_j(\zeta)|^s) \to 0 ~\mbox{as}~ j \to \infty
		\end{eqnarray*}
		implies that $|Dg(\zeta)|\equiv 0.$ That is,
		$$\frac{\partial g(\zeta)}{\partial z_i}= 0, ~i= 1, 2, \ldots, n.$$ which shows that $g$ is constant, a contradiction. $\Box$\\
	
\bibliographystyle{amsplain}



\end{document}